\documentclass[12pt,oneside]{amsart}
\usepackage{geometry}                
\usepackage{graphicx}
\usepackage{amssymb}
\usepackage{epstopdf}

\usepackage{amsmath,amsthm,amscd,amssymb}
\usepackage{latexsym}
\usepackage[colorlinks,citecolor=red,pagebackref,hypertexnames=false]{hyperref}
\numberwithin{equation}{section}
\theoremstyle{plain}
\newtheorem{theorem}{Theorem}[section]
\newtheorem{lemma}[theorem]{Lemma}
\newtheorem{corollary}[theorem]{Corollary}

\theoremstyle{definition}
\newtheorem{definition}[theorem]{Definition}

\theoremstyle{remark}
\newtheorem{remark}[theorem]{Remark}

\newtheorem{case[theorem]}{Case}

\date{September 6, 2014}      
\author{M. Bennett, A. Iosevich and K. Taylor}
\address{Department of Mathematics, University of Rochester, Rochester, NY}
\email{bennett@math.rochester.edu}
\address{Department of Mathematics, University of Rochester, Rochester, NY}
\email{iosevich@math.rochester.edu}
\address{Institute for Mathematics and its Applications}
\email{kltaylor@ima.umn.edu} 
\thanks{This work was partially supported by the NSF Grant DMS10-45404.}
\title{\parbox{14cm}{\centering{Finite chains inside thin subsets of ${\mathbb R}^ \MakeLowercase{d}$}}}
\begin{document}
\maketitle
\begin{abstract} In a recent paper, Chan, \L aba, and Pramanik investigated geometric configurations inside thin subsets of the Euclidean set possessing measures with Fourier decay properties. In this paper we ask which configurations can be found inside thin sets of a given Hausdorff dimension without any additional assumptions on the structure. We prove that if the Hausdorff dimension of $E \subset {\Bbb R}^d$, $d \ge 2$, is greater than $\frac{d+1}{2}$, then for each $k \in \mathbb{Z}^+$ there exists a non-empty interval $I$ such that, given any sequence $\{t_1, t_2, \dots, t_k; t_j \in I\}$, there exists a sequence of distinct points ${\{x^j\}}_{j=1}^{k+1}$, such that $x^j \in E$ and $|x^{i+1}-x^i|=t_j$, $1 \leq i \leq k$. In other words, $E$ contains vertices of a chain of arbitrary length with prescribed gaps.
\end{abstract}  

\maketitle

\section{Introduction} The problem of determining which geometric configurations one can find inside various subsets of Euclidean space is a classical subject matter. The basic problem is to understand how large a subset of Euclidean space must be to be sure that it contains the vertices of a congruent and possibly scaled copy of a given polyhedron or another geometric shape. In the case of a finite set, ``large" refers to the number of points, while in infinite sets, it refers to the Hausdorff dimension or Lebesgue density. The resulting class of problems has been attacked by a variety of authors using combinatorial, number theoretic, ergodic, and Fourier analytic techniques, creating a rich set of ideas and interactions. 

We begin with a comprehensive result due to Tamar Ziegler, \cite{Z06} which generalizes an earlier result due to Furstenberg, Katznelson and Weiss \cite{FKW90}. See also \cite{B86}. 

\begin{theorem} \label{z06} [Ziegler] Let $E \subset {\Bbb R}^d$, of positive upper Lebesgue density in the sense that 
$$ \limsup_{R \to \infty} \frac{{\mathcal L}^d \{E \cap {[-R,R]}^d \}}{{(2R)}^d}>0, $$ where ${\mathcal L}^d$ denotes the $d$-dimensional Lebesgue measure. Let $E_{\delta}$ denote the $\delta$-neighborhood of $E$. Let $V=\{ {\bf 0}, v^1, v^2, \dots, v^{k-1}\} \subset {\Bbb R}^d$, where $k \ge 2$ is a positive integer. Then there exists $l_0>0$ such that for any $l>l_0$ and any $\delta>0$ there exists $\{x^1, \dots, x^k\} \subset E_{\delta}$ congruent to $lV=\{ {\bf 0}, lv^1, \dots, 
lv^{k-1}\}$. \end{theorem} 

In particular, this result shows that we can recover every simplex similarity type and sufficiently large scaling inside a subset of ${\Bbb R}^d$ of positive upper Lebesgue density. It is reasonable to wonder whether the assumptions of Theorem \ref{z06} can be weakened, but the following result due to Maga \cite{Mag10} shows that conclusion may fail even if we replace the upper Lebesgue density condition with the assumption that the set is of dimension $d$. 

\begin{theorem} \label{mag10} [Maga] For any $d \ge 2$ there exists a full dimensional compact set $A \subset {\Bbb R}^d$ such that $A $ does not contain the vertices of any parallelogram. If $d=2$, then given any triple of points $x^1,x^2,x^3$, $x^j \in A$, there exists a full dimensional compact set $A \subset {\Bbb R}^2$ such that $A$ does not contain the vertices of any triangle similar to $\bigtriangleup x^1x^2x^3$. \end{theorem} 

In view of Maga's result, it is reasonable to ask whether interesting point configurations can be found inside thin sets under additional structural hypotheses. This question was recently addressed by Chan, \L aba, and Pramanik in \cite{CLP14}.  Before stating their result, we provide two relevant definitions. 

\begin{definition}\label{clpConfiguration} 
Fix integers $n\geq 2$, $p\geq 3$, and $m= n\lceil \frac{p+1}{2} \rceil$.  Suppose $B_1, \dots, B_p$ are $n \times (m-n)$ matrices.

(a)  We say that $E$ contains a $p-$point $\mathcal{B}-$configuration if there exists vectors
$z\in \mathbb{R}^n $ and $w\in \mathbb{R}^{m-n}\backslash \vec{0}$ such that $$\{z + B_j w \}_{j=1}^p \subset E.$$ 

(b) Moreover, given any finite collection of subspaces $V_1,\dots, V_q \subset \mathbb{R}^{m-n}$ with $dim(V_i) < m-n$, we say that $E$ contains a non-trivial $p-$point $\mathcal{B}-$configuration with respect to $(V_1,\dots, V_q)$ if there exists vectors
$z\in \mathbb{R}^n$ and $w\in \mathbb{R}^{m-n}\backslash \cup_{i=1}^{q}V_i$ such that $$\{z + B_j w \}_{j=1}^p \subset E.$$ \\
\end{definition}

\begin{definition}\label{clpRank}
Fix integers $n\geq 2$, $p\geq 3$, and $m= n\lceil \frac{p+1}{2} \rceil$. We say that a set of $n\times (m-n)$ matrices $\{ B_1, \dots, B_p\}$ is non-degenerate if  

\[rank \left( \begin{array}{c}
B_{i_2}-B_{i_1}\\
\vdots\\
B_{i_{m/n}}-  B_{i_1}\\
\end{array} \right)=m-n \]

for any distinct indices $i_1,\dots,i_{ m/n} \in \{1,\dots,p\}$.
\end{definition}

\vskip.125in 

\begin{theorem} \label{clp} [Chan, \L aba, and Pramanik] 
Fix integers $n\geq 2$, $p\geq 3$, and $m= n\lceil \frac{p+1}{2} \rceil$.  Let $\{B_1, \dots, B_p\}$ be a collection of $n \times (m-n)$ non-degenerate matrices in the sense of Definition \ref{clpRank}.  Then for any constant $C$, there exists a positive number $\epsilon_0 = \epsilon_0(C,n,p,B_1,\dots,B_p) <<1$ with the following property:  Suppose the set $E \subset \mathbb{R}^n$ with $\left|E \right|=0$ supports a positive, finite, Radon measure $\mu$ with two conditions:
(a) (ball condition)  $sup_{\stackrel{x\in E}{ 0<r<1}} \frac{\mu(B(x,r)}{r^{\alpha}} \le C$ if $n-\epsilon_0 <\alpha < n$,
(b) (Fourier decay) $sup_{\xi \in \mathbb{R}^n } |\widehat{\mu}(\xi)| (1+ |\xi|)^{\beta/2} \le C.$

Then

\vskip.125in 

(i) $E$ contains a $p-$point $\mathcal{B}-$configuration in the sense of Definition \ref{clpConfiguration} (a).

\vskip.125in 

(ii) Moreover, for any finite collection of subspaces $V_1, \dots, V_q \subset \mathbb{R}^{m-n}$ with 

$dim(V_i) < m-n$, $E$ contains a non-trivial $p-$point $\mathcal{B}-$configuration with respect 

to $(V_1, \dots, V_q)$ in the sense of Definition \ref{clpConfiguration} (b). 

\end{theorem} 

\vskip.125in 

One can check that the Chan-\L aba-Pramanik result covers some geometric configurations but not others. For example, their non-degeneracy condition allows them to consider triangles in the plane, but not simplexes in ${\Bbb R}^3$ where three faces meet at one of the vertices at right angles, forming a three-dimensional corner. Most relevant to this paper is the fact that the conditions under which Theorem \ref{clp} holds are satisfied for chains (see Definition \ref{chaindefinition} below), but the conclusion requires decay properties for the Fourier transform of a measure supported on the underlying set. We shall see that in the case of chains, such an assumption is not needed and the existence of a wide variety of chains can be established under an explicit dimensional condition alone. 

\subsection{Focus of this article} 

In this paper we establish that a set of sufficiently large Hausdorff dimension, {\it with no additional assumptions}, contains an arbitrarily long chain with vertices in the set and preassigned admissible gaps.  

\vskip.125in 

\begin{definition} \label{chaindefinition} (See Figure 1 above) A $k$-chain in $E \subset {\Bbb R}^d$ with gaps ${\{t_i\}}_{i=1}^k$ is a sequence 
$$\{x^1,x^2, \dots, x^{k+1}: x^j \in E; \ |x^{i+1}-x^i|=t_i; \ 1 \leq i \leq k\}.$$ 

\vskip.125in 

We say that the chain is {\it non-degenerate} if all the $x^j$s are distinct. 
\end{definition} 

\begin{figure}
\label{chainfigure}
\centering
\includegraphics[scale=.5]{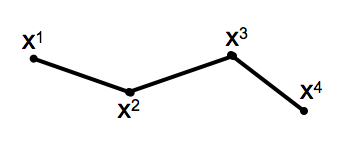}
\caption{A 3-chain}
\end{figure}

\vskip.125in

Our main result is the following. 

\begin{theorem} \label{main} Suppose that the Hausdorff dimension of a compact set $E \subset {\Bbb R}^d$, $d \ge 2$, is greater than $\frac{d+1}{2}$. Then for any $k \ge 1$, there exists an open interval $\tilde{I}$ such that for any  ${\{t_i\}}_{i=1}^k \subset \tilde{I}$ there exists a non-degenerate $k$-chain in $E$ with gaps  ${\{t_i\}}_{i=1}^k$. \end{theorem} 

\vskip.125in 

In the course of establishing Theorem \ref{main} we shall prove the following result which is interesting in its own right and has a number of consequences for Falconer type problems. See \cite{Fal86}, \cite{Erd05} and \cite{W99} for the background and the latest results pertaining to Falconer distance problem. 

\begin{theorem} \label{almostmain} Suppose that $\mu$ is a compactly supported non-negative Borel measure such that 
\begin{equation} \label{adupper} \mu(B(x,r)) \leq Cr^{s_{\mu}}, \end{equation} where $B(x,r)$ is the ball of radius $r>0$ centered at $x \in {\Bbb R}^d$, for some $s_{\mu}\in(\frac{d+1}{2}, d]$. Then for any $t_1, \dots, t_k>0$ and $\epsilon>0$,
\begin{equation} \label{cbabove} \mu \times \mu \times \dots \times \mu \{(x^1,x^2, \dots, x^{k+1}): t_i -\epsilon \leq |x^{i+1}-x^i| \leq t_i+\epsilon; \ i=1,2, \dots, k \} \leq C\epsilon^k. \end{equation} \end{theorem} 

\vskip.125in 

\begin{corollary} \label{chainfalconer} Given a compact set $E \subset {\Bbb R}^d$, $d \ge 2$, $k\geq 1$, define 
$$ \Delta_k(E)=\left\{|x^1-x^2|, |x^2-x^3|, \dots, |x^k-x^{k+1}|: x^j \in E \right\}.$$ 

Suppose that the Hausdorff dimension of $E$ is greater than $\frac{d+1}{2}$. Then 
$${\mathcal L}^k(\Delta_k(E))>0.$$
\end{corollary} 

\vskip.125in 


\vskip.125in 

\begin{remark}
Suppose that $E \subset {\Bbb R}^d$ has Hausdorff dimension $s>\frac{d+1}{2}$ and is \textbf{Ahlfors-David regular}, i.e. there exists $C>0$ such that for every $x \in E$,
$$ C^{-1}r^s \leq \mu(B(x,r)) \leq Cr^s,$$ 
(where $\mu$ is the restriction of the $s$-dimensional Hausdorff measure to $E$). Then using the techniques in \cite{EIT11} along with Theorem \ref{almostmain}, one can show that for any sequence of positive real numbers $t_1, t_2, \dots, t_k$, the upper Minkowski dimension of 
$$ \{(x^1, x^2, \dots, x^{k+1}) \in E^{k+1}: |x^{j+1}-x^j|=t_j; \ 1 \leq j \leq k\}$$
does not exceed $(k+1)dim_{{\mathcal H}}(E)-k$.
\end{remark}

\vskip.125in 

\subsection{Acknowledgements} The authors wish to thank Shannon Iosevich for her help with the diagrams used in this paper. The authors also wish to thank Fedja Nazarov and Jonathan Pakianathan for helpful discussions related to the subject matter of this article. 

\vskip.25in 

\section{Proof of Theorem \ref{main} and Theorem \ref{almostmain}} 

\vskip.125in 
The strategy for this section is as follows:

\vspace{.25 in}
We begin by dividing both sides of equation \eqref{cbabove} by $\epsilon^k$.  The left side becomes
\begin{equation}\label{density}
\epsilon^{-k}\mu \times \dots \times \mu \{(x^1,\dots, x^{k+1}): t_i-\epsilon \leq |x^{i+1}-x^i| \leq t_i+\epsilon; \ i=1,2, \dots, k \},
\end{equation}
which can be interpreted as the density of $\epsilon$-approximate chains in $E \times \ldots \times E$.

\vspace{.125 in}
\noindent
Theorem \ref{almostmain} gives an upper bound on this expression that is independent of $\epsilon$. This is accomplished using an inductive argument on the chain length coupled with repeated application of an earlier result from \cite{ISTU14} in which the authors establish $L^2(\mu)$ mapping properties of certain convolution operators. 
This upper bound is important in the final section where we define a measure on the set of chains.  

Next, we acquire a lower bound on \eqref{density}.  This result was already established in the case that $k=1$ in \cite{IMT12} where the authors show that the density of $\epsilon$-approximate $1$-chains with gap size $t$ is bounded below independent of $\epsilon$ for all $t$ in a non-empty open interval, $I$. 
Using a pigeon-holing argument, we extend the result in \cite{IMT12} to obtain a lower bound on \eqref{density} in the case that every gap is of equal size, $t$, for some $t\in I$.  
To obtain a lower bound on chains with variable gap size, we show that the density of $\epsilon$-approximate $k$-chains is continuous as a function of gap sizes. Furthermore, we use the lower bound on chains with constant gaps to prove that this continuous function is not identically zero.  We conclude that the density of $\epsilon$-approximate $k$-chains is bounded below independent of $\epsilon$ and independent of the gap sizes, as long as all gap sizes fall within some interval $\tilde{I}$ around $t$.
\vskip.125in

In the final section, we address the issue of non-degeneracy.  To this end, we reinterpret the density of $\epsilon$-approximate $k$-chains as a measure supported in $E^{k+1}$, and show that it converges to a new measure, $\Lambda_{\vec{t}}^k$, as $\epsilon \downarrow 0$. This new measure is shown to be supported on ``exact" $k$-chains ($\epsilon = 0$) with admissible gaps.  We next show that the measure of the set of degenerate chains is $0$, and we conclude that the mass of $\Lambda_{\vec{t}}^k$ is contained in non-degenerate $k-$chains. 

\vspace{.25 in}

We shall repeatedly use the following result due to Iosevich, Sawyer, Taylor, and Uriarte-Tuero \cite{ISTU14}. 

\begin{theorem} \label{maintool} Let $T_{\lambda}f(x)=\lambda*(f\mu)(x)$, where $\lambda, \mu$ are compactly supported non-negative Borel measures on ${\Bbb R}^d$. Suppose that $\mu$ satisfies \eqref{adupper} and for some $\alpha>0$
$$ |\widehat{\lambda}(\xi)| \leq C{|\xi|}^{-\alpha}.$$ 
Suppose that $\nu$ is a compactly supported Borel measure supported on ${\Bbb R}^d$ satisfying \eqref{adupper} with $s_{\mu}$ replaced by $s_{\nu}$ and suppose that $\alpha>d-s$, where $s=\frac{s_{\mu}+s_{\nu}}{2}$. Then 
$$ {||T_{\lambda}f||}_{L^2(\nu)}  \leq c {||f||}_{L^2(\mu)}.$$
\end{theorem} 

\vskip.125in

In this article, we will use Theorem \eqref{maintool} with $\lambda=\sigma,$ the surface measure on a $(d-1)$-dimensional sphere in $\mathbb{R}^d$.  It is known, see  \cite{St93}, that $$\widehat{\sigma}(\xi) =O(|\xi|^{-(d-1)/2}).$$

\vskip.125in

Since the proof of Theorem \ref{maintool} is short, we give the argument below for the sake of keeping the presentation as self-contained as possible. It is enough to show that 
$$ \langle T_{\lambda^{\epsilon}}f, g\nu \rangle  \leq C{||f||}_{L^2(\mu)} \cdot {||g||}_{L^2(\nu)}.$$ 

The left hand side equals 
$$ \int \hat{\lambda}^{\epsilon}(\xi) \widehat{f\mu}(\xi) \widehat{g\nu}(\xi) d\xi.$$ By the assumptions of Theorem \ref{maintool}, the modulus of this quantity is bounded by 
$$ C \int {|\xi|}^{-\alpha} |\widehat{f\mu}(\xi)| |\widehat{g\nu}(\xi)| d\xi,$$ and applying Cauchy-Schwarz bounds this quantity by 
\begin{equation} \label{2square} C {\left( \int {|\widehat{f\mu}(\xi)|}^2 {|\xi|}^{-\alpha_{\mu}} d\xi \right)}^{\frac{1}{2}} \cdot 
{\left( \int {|\widehat{g \nu}(\xi)|}^2 {|\xi|}^{-\alpha_{\nu}} d\xi \right)}^{\frac{1}{2}} \end{equation} for any $\alpha_{\mu}, \alpha_{\nu}>0$ such that $\alpha=\frac{\alpha_{\mu}+\alpha_{\nu}}{2}$. 

\vskip.125in 

By Lemma (\ref{fenergy}) below, the quantity (\ref{2square}) is bounded by $C {||f||}_{L^2(\mu)} \cdot {||g||}_{L^2(\nu)}$ after choosing, as we may, $\alpha_{\mu}>d-s_{\mu}$ and $\alpha_{\nu}>d-s_{\nu}$. This completes the proof of Theorem \ref{maintool}. 

\vskip.25in




\vskip.25in 

\subsection{Proof of Theorem \ref{almostmain} and Corollary \ref{chainfalconer}}
Let $\epsilon>0$.  Divide both sides of (\ref{cbabove}) by $\epsilon^k$, and note that it suffices to establish the estimate 
\begin{equation}\label{up} C_k^{\epsilon}(\mu)=\int \left( \prod_{i=1}^k \sigma_{t_i}^{\epsilon}(x^{i+1}-x^i) d\mu(x^i) \right) d\mu(x^{k+1}) \le c^k,\end{equation}
where $c$ is independent of $\epsilon$, and $t_1, \ldots,t_ k>0$. Here 
$\sigma_r^{\epsilon}(x)=\sigma_r*\rho_{\epsilon}(x)$, with $\sigma_r$ the Lebesgue measure on the sphere of radius $r$, $\rho$ a smooth cut-off function with $\int \rho=1$ and $\rho_{\epsilon}(x)=\epsilon^{-d} \rho \left( \frac{x}{\epsilon} \right)$. Assume in addition that $\rho$ is non-negative and that $\rho(x) =\rho(-x)$.
\vskip.125in

\vskip.125in

Let $\sigma$ denote the Lebesgue measure on the $(d-1)-$dimensional sphere in $\mathbb{R}^d$.  Set $T_j^\epsilon = T_{\sigma_{t_{j}}}^\epsilon$, where $T_{\sigma_{t_j}}^{\epsilon}f(x) = \sigma_{t_j}*(f\mu)(x)$ was introduced in Theorem \ref{maintool}.
Define
	\begin{equation}\label{fk}f_k^{\epsilon}(x) = T_k^\epsilon \circ \cdots \circ T_1^{\epsilon}(1)(x),\end{equation} 
and $$f_0^{\epsilon}(x)=1.$$
It is important to note that $f_k(x)$ depends implicitly on  the choices of $t_1, \dots, t_k>0$, and this choice will be made explicit throughout.  
 \vskip.125in

Observe that 
\begin{equation}\label{beauty}f_{k+1}^{\epsilon}=T_{k+1}^{\epsilon}f_k^{\epsilon}.\end{equation}
\vskip.125in

Re-writing the left-hand-side of \eqref{up}, it suffices to show
\begin{equation}\label{upup}C_k^{\epsilon}(\mu) = \int f_k^{\epsilon}(x) d\mu(x)\le c^k.\end{equation}

Using Cauchy-Schwarz (and keeping in mind that $\int d\mu(x) = 1$), we bound the left-hand-side of \eqref{upup} by
\begin{equation}C_k^{\epsilon}(\mu) = \int f_k^{\epsilon}(x) d\mu(x)\le   \|f_k^{\epsilon}\|_{L^2(\mu)}        .\end{equation}

We now use induction on $k$ to show that 
\begin{equation}\label{mini}\|f_k^{\epsilon}\|_{L^2(\mu)} \le c^{k},\end{equation} where $c$ is the constant obtained in Theorem \ref{maintool}.
For the base case, $k=0$, we have $\|f_0^{\epsilon}\|_{L^2(\mu)}=  \int d\mu(x)  =1.$
Next, we assume inductively  that $\|f_k^{\epsilon}\|_{L^2(\mu)} \le c^{k}$.  

We now show that, for any $t_{k+1}>0$, $$\|f_{k+1}^{\epsilon}\|_{L^2(\mu)}\le c^{k+1}.$$ 

First, use \eqref{beauty} to write 
$$\|f_{k+1}^{\epsilon}\|_{L^2(\mu)}=||T_{k+1}^{\epsilon}f_k^\epsilon||_{L^2(\mu)}.$$
Next, use Theorem \ref{maintool} with $\lambda=\sigma$, the Lebesgue measure on the sphere,  and $\alpha = \frac{d-1}{2}$ (see the comment immediately following Theorem 2.1 to justify this choice of $\alpha$) to show that 
$$||T_{k+1}^{\epsilon}f_k^\epsilon||_{L^2(\mu)} \le c\|f_{k}^{\epsilon}\|_{L^2(\mu)}$$
whenever $s_{\mu}>d-\alpha =\frac{d+1}{2}$.  

We complete the proof by applying the inductive hypothesis.  This completes the verification of \eqref{mini}. 
\\

We now recover Corollary \ref{chainfalconer}. Let $s_{\mu} \in \left(\frac{d+1}{2}, dim(E)\right)$, and choose a probability measure, $\mu$, with support contained in $E$  which satisfies \eqref{adupper}; the existence of such a measure is provided by Frostman's lemma (see \cite{Falc86}, \cite{W04} or \cite{M95}).

Cover $\Delta_k(E)$ with cubes of the form 
$$ \bigcup_i \prod_{j=1}^d (t_{ij}, t_{ij}+\epsilon_i),$$ where $\prod$ denotes the Cartesian product. We have 

\begin{align*} \label{covering} 
1&=\mu\times \cdots \times \mu(E^{k+1})\\
& \leq \sum_i 
\mu \times \cdots \times \mu \{(x^1, \dots, x^{k+1}): t_{ij} -\epsilon\leq |x^{j+1}-x^j| \leq t_{ij}+\epsilon_i; \ 1 \leq j \leq k \}. \end{align*}
By Theorem \ref{almostmain},  the expression above is bounded by 
\begin{equation} \label{ksum} C \sum_i \epsilon_i^k. \end{equation} and we conclude that (\ref{ksum}) is bounded from below by $\frac{1}{C}>0$. It follows that $\Delta_k(E)$ cannot have measure $0$ and the proof of Corollary \ref{chainfalconer} is complete. 

\vskip.25in 

We now continue with the proof of Theorem \ref{main}. 

\subsection{Lower bound on $C^\epsilon_k(\mu)$}
 Let $s_{\mu} \in \left(\frac{d+1}{2}, dim(E)\right)$, and choose a probability measure, $\mu$, with support contained in $E$  which satisfies \eqref{adupper}.
\vskip.125in

We now establish the existence of a non-empty open interval $\tilde{I}$ such that  
\begin{equation} \label{lowerboundest} \liminf_{\epsilon \to 0} C^{\epsilon}_k(\mu) >0 \end{equation}
where each $t_i$ belongs to $\tilde{I}$, $C^\epsilon_k(\mu)$ is as in \eqref{up}.

Note that this positive lower bound alone establishes the existence of vertices $x^1, \dots, x^{k+1} \in E$ so that $|x^{i+1}-x^{i}|=t_i$ for each $i\in \{1,\dots, k\}$ (this follows, for instance, by Cantor's intersection theorem and the compactness of the set $E$).  Extra effort is made in the next section in order to guarantee that we may take $x^1, \dots, x^{k+1}$ distinct. 

We first prove estimate \eqref{lowerboundest} in the case that all gaps are equal. This is accomplished using a pigeon holing argument on chains of length one.  
We then provide a continuity argument to show that the estimate holds for variable gap values, $t_i$, belonging to a non-empty open interval $\tilde{I}$.
The second argument relies on the first precisely at the point when we show that the said continuous function is not identically equal to zero. 
\vskip.125in

\noindent
\textbf{Lower bound for constant gaps:}
\\
 The proof of estimate \eqref{lowerboundest}  in the case when $k=1$ was already established in \cite{IMT12} provided that $\mu$ to satisfies the ball condition in \eqref{adupper} with $\frac{d+1}{2}<s_{\mu}<dim_{{\mathcal H}}(E)$. The existence of such measures is established by Frostman's lemma (see e.g. \cite{Falc86}, \cite{W04} or \cite{M95} ).  

More specifically, it is demonstrated in \cite{IMT12} that there exists $c(1)>0$, $\epsilon_0>0$, and a non-empty open interval $I \subset (0, diameter(E))$ so that if $t\in I$ and $0<\epsilon < \epsilon_0$, then  
$$C_1^{\epsilon}= \int \sigma_t^{\epsilon}*\mu(x) d\mu(x)> 2c(1).$$ 

To establish estimate \eqref{lowerboundest} for longer chains, we rely on the following lemmas.  
\begin{lemma}\label{babylower}
Set $$G_{t,\epsilon}(1) = \{x\in E: \sigma_t^{\epsilon}*\mu(x) >c(1)\}.$$
There exists $m(1) \in \mathbb{Z}^+$ so that if $t\in I$ and $0<\epsilon <\epsilon_0$, then 
$$\mu( G_{t,\epsilon}(1) )    \geq    2^{-2m(1)}.$$
\end{lemma}

\vskip.125in 


\begin{lemma}\label{mamalower}
 Set $$G_{t,\epsilon}(j+1) = \{x\in E: \sigma_t^{\epsilon}*\mu|_j(x) >c(j+1)\},$$
where $j \in \{1,\cdots, (k-1)\}$, $\mu|_j(x)$ denotes restriction of the measure $\mu$ to the set ${G_{t, \epsilon }(j)}$, and 
$$c(j+1) =\frac{1}{2}c(j)\mu(G_{t,\epsilon}(j)).$$ 
Then there exists $m(j+1) \in \mathbb{Z}^+$ so that if $t\in I$ and $0<\epsilon <\epsilon_0$, then 
$$\mu( G_{t,\epsilon}(j+1) )> 2^{-2m(j+1)}.$$

\end{lemma}
\vskip.125in

We postpone the proof of Lemmas \ref{babylower} and \ref{mamalower} momentarily, and we apply these lemmas to obtain a lower bound on $C_{k}^{\epsilon}(\mu)$.  
\\

We write 
$$C_{k}^{\epsilon}(\mu)
= \int f_k^{\epsilon}(x)d\mu(x),$$ 
where $f_k^{\epsilon}$ was introduced in \eqref{fk} and here $t_1=\cdots =t_k =t$.  
\vskip.125in 

Now
\begin{align*}
C_{k}^{\epsilon}(\mu)
=& \int f_k^{\epsilon}(x)d\mu(x)\\
=& \iint \sigma_{t}^{\epsilon}(x-y) f_{k-1}(y) d\mu(y) d\mu(x)\\
\end{align*}

 Integrating in $x$ and restricting the variable $y$ to the set $G_{t,\epsilon}(1)$, we write
\begin{align*}
 C_{k}^{\epsilon}(\mu) 
 \geq & \int_{G_{t, \epsilon}(1)} \sigma_t^{\epsilon}*\mu(y) f_{k-1}(y) d\mu(y)
\\
 \geq & \,  c(1)  \int_{G_{t, \epsilon}(1)}  f_{k-1}(y) d\mu(y)
\\
 =& \,  c(1)  \int  f_{k-1}(y) d\mu_1(y).
\end{align*}

To achieve a lower bound, we iterate this process. For each $j\in \{2, \dots, k-1\}$ we have:
\begin{align*}
& \int f_{k-j}^{\epsilon}(x)d\mu_j(x)
\\
= & \iint \sigma_{t}^{\epsilon}(x-y) f_{k-j-1}(y) d\mu(y) d\mu_j(x)
\\
 \geq & \int_{G_{t, \epsilon}(j+1)} \sigma_t^{\epsilon}*\mu_j(y) f_{k-j-1}(y) d\mu(y)
\\
 \geq & \,  c(j+1)  \int_{G_{t, \epsilon}(j+1)}  f_{k-j-1}(y) d\mu(y)
\\
 =& \,  c(j+1)  \int f_{k-j-1}(y) d\mu_{j+1}(y)
\\
\end{align*}
It follows that 
\begin{align*}
 C_{k}^{\epsilon}(\mu) &\geq \left(\prod_{j=1}^{k-1}c(i)\right)   \int  \int \sigma_t^{\epsilon}(x-y) d\mu_{k-1}(y)  d\mu(x)\\
&\geq \left(\prod_{j=1}^{k}c(i)\right)  \mu(G_{t,\epsilon}(k)),\\
\end{align*}
and we are done in light of Lemma \ref{mamalower}.

Given Lemmas \ref{babylower} and \ref{mamalower}, we have shown that for all $t\in I$ and for all $0<\epsilon<\epsilon_0$, we have
\begin{equation}\label{upperconstant}\liminf_{\epsilon \to 0} C_k^{\epsilon}(\mu)> 0,\end{equation}
where all gap lengths, $t_1, \dots, t_k$ constantly equal to $t$.
This concludes the proof of estimate \eqref{lowerboundest} in the case of constant gaps.  
\vskip.125in

\noindent
We now proceed to the proofs of Lemmas \ref{babylower} and \ref{mamalower}.

\begin{proof}(Lemma \ref{babylower})\\
We write 
\begin{align*}
 2c(1)  &<   \int  \sigma_t^{\epsilon}*\mu(x)   d\mu(x)   \\
& \le \left(\int_{(G_{t,\epsilon}(1))^c} \sigma_t^{\epsilon}*\mu(x)   d\mu(x) \right)  + \left( \int_{G_{t,\epsilon}(1)} \sigma_t^{\epsilon}*\mu(x)   d\mu(x) \right) 
\\
&=\mathcal{I} + \mathcal{II}\\
\end{align*}
where $A^c$ denotes the compliment of a set $A\subset E.$

\vspace{.125 in}
\noindent
We first observe that $$\mathcal{I} \le c(1).$$  

Next, we estimate $\mathcal{II}$.  Let $m\in \mathbb{Z}^+$, and write 
$$G_{t,\epsilon}(1) =  \{x\in E: c(1) < \sigma_t^{\epsilon}*\mu(x) \le 2^m\} \cup  \{x\in E: 2^m \le \sigma_t^{\epsilon}*\mu(x) \}.$$
Then
\begin{align*}
\mathcal{II} &= \left( \int_{\{x\in E: c(1) < \sigma_t^{\epsilon}*\mu(x) \le 2^m\}} \sigma_t^{\epsilon}*\mu(x)   d\mu(x) \right) +  \left( \int_{\{x\in E: 2^m \le  \sigma_t^{\epsilon}*\mu(x) \}} \sigma_t^{\epsilon}*\mu(x)   d\mu(x) \right) 
\\
& \le 2^m \mu( G_{t,\epsilon}(1))  +\left( \sum_{l=m} 2^{l+1}\cdot \mu ( \{x\in E: 2^l \le \sigma_t^{\epsilon}*\mu(x) \le 2^{l+1}\}) \right) .
\\
\end{align*}

We use Theorem \ref{maintool} to estimate $$\mu ( \{x\in E: 2^l \le \sigma_t^{\epsilon}*\mu(x) \le 2^{l+1}\})  \le c_d \cdot 2^{-2l},$$ where the constant $c_d$ depends only on the ambient dimension $d$.  
Now, 
\begin{align*}
\mathcal{II}
&\le 2^m \mu( G_{t,\epsilon}(1)) +  \left(2 c_d \cdot \sum_{l=m} 2^{l}\cdot 2^{-2l}\right)  \\
& \lesssim 2^m \mu( G_{t,\epsilon}(1))  +  2^{-m}. 
\\
\end{align*}
It follows that 
$$2c(1) \le \mathcal{I} + \mathcal{II} \lesssim   c(1) + 2^m \mu( G_{t,\epsilon}(1))  +  2^{-m}. $$
Taking $m \in \mathbb{Z}^+$ large enough, we conclude that $$\mu( G_{t,\epsilon}(1)) \geq 2^{-2m}.$$
\end{proof}

\begin{proof}(Lemma \ref{mamalower})\\
We prove the Lemma by induction on $j$.  The base case, $j=1$, was established in Lemma \ref{babylower}.  Next, assume that there exists $ m(j) \in \mathbb{Z}^+$ such that 
$$2^{-m(j)} < \mu(G_{t,\epsilon}(j) )$$ for all $0<\epsilon< \epsilon_0$ and $t\in I$.  \\

By the definition of $G_{t, \epsilon}(j)$,
\begin{align*}    c(j)\mu(G_{t,\epsilon}(j)) &< \int_{G_{t, \epsilon}(j)} \sigma^{\epsilon}_t*\mu|_{G_{t,\epsilon}(j-1)} (x) d\mu(x)  .
\\
\end{align*}
Set $c(j+1) = \frac{1}{2}c(j)\mu(G_{t,\epsilon}(j))$.  By assumption, $2c(j+1)= c(j)\mu(G_{t,\epsilon}(j) \geq c(j)2^{-m(j)}$, and in particular this quantity is positive.
Next, we obtain a bound from above:  
\begin{align*}  \int_{G_{t, \epsilon}(j)} \sigma^{\epsilon}_t*\mu|_{G_{t,\epsilon}(j-1)} (x) d\mu(x) 
& \le   \int_{G_{t, \epsilon}(j)} \sigma^{\epsilon}_t*\mu (x) d\mu(x) 
\\
&= \int \sigma^{\epsilon}_t*\mu|_{j} (x) d\mu(x)
\\
&= \left(\int_{(G_{t, \epsilon}(j+1))^c} \sigma^{\epsilon}_t*\mu|_{j} (x) d\mu(x)\right) + \left(\int_{G_{t, \epsilon}(j+1)} \sigma^{\epsilon}_t*\mu|_{j} (x) d\mu(x)\right)
\\
&=\mathcal{I} + \mathcal{II}.
\\
\end{align*}
First we observe that $$\mathcal{I} \le c(j+1).$$  

Next, we estimate $\mathcal{II}.$ Let $m\in \mathbb{Z}^+$, and write 
$$G_{t,\epsilon}(j+1) =  \{x\in E: c(j+1) < \sigma_t^{\epsilon}*\mu|_{j}(x) \le 2^m\} \cup  \{x\in E: 2^m \le \sigma_t^{\epsilon}*\mu|_{j}(x) \}.$$
Then
\begin{align*}
\mathcal{II} &= \left( \int_{\{x\in E: c(j+1) < \sigma_t^{\epsilon}*\mu|_{j}(x) \le 2^m\}} \sigma_t^{\epsilon}*\mu|_{j}(x)   d\mu(x) \right) +  \left( \int_{\{x\in E: 2^m \le  \sigma_t^{\epsilon}*\mu(x) \}} \sigma_t^{\epsilon}*\mu|_{j}(x)   d\mu(x) \right) 
\\
& \le 2^m \cdot \mu( G_{t,\epsilon}(j+1)) +\left( \sum_{l=m} 2^{l+1}\cdot \mu ( \{x\in E: 2^l \le \sigma_t^{\epsilon}*\mu|_{j}(x) \le 2^{l+1}\}) \right) .
\\
\end{align*}

We use Theorem \ref{maintool} to estimate $$\mu ( \{x\in E: 2^l \le \sigma_t^{\epsilon}*\mu|_{j}(x) \le 2^{l+1}\})  \le c_d \cdot 2^{-2l},$$ where the constant $c_d$ depends only on the ambient dimension $d$ and the choice of the measure $\mu$.  
Now, 
\begin{align*}
\mathcal{II} 
&\le 2^m \mu( G_{t,\epsilon}(j+1)) +  \left(2 c_d \cdot \sum_{l=m} 2^{l}\cdot 2^{-2l}\right)  \\
& \lesssim 2^m \mu( G_{t,\epsilon}(j+1))  +  2^{-m}. 
\\
\end{align*}
It follows that 
$$2c(j+1) \le \mathcal{I} +\mathcal{II} \lesssim   c(j+1) + 2^m \mu( G_{t,\epsilon}(j+1))  +  2^{-m}. $$
Taking $m \in \mathbb{Z}^+$ large enough, we conclude that $$\mu( G_{t,\epsilon}(j+1)) \geq 2^{-2m}.$$
\end{proof}
\vskip.125in

\noindent
\textbf{Lower bound for variable gaps}
\smallskip

%
%
We now verify \eqref{lowerboundest} in the case of variable gap lengths.  In more detail, we show that, for all $k \in \mathbb{Z}^+$ and for values of $t_i$ in a non-empty open interval $\tilde{I}$, we have
\begin{equation}\label{repeat}\liminf_{\epsilon \to 0} \int f_k^\epsilon (x)d\mu(x) > 0,\end{equation}
where $f_k^{\epsilon}$ is defined in \eqref{fk} with $0<t_1,\dots, t_k \in \tilde{I}$.

\vskip.125in 
The following lemma captures the strategy of proof, and establishes \eqref{repeat}.

\begin{lemma} \label{epsLimit}
\begin{equation}\label{rewrite}C^{\epsilon}_k(\mu)=\int f_{k}^\epsilon (x)d\mu(x) =M_k(t_1, \dots, t_k) - \sum_{j=1}^{k} R_{k,j}^{\epsilon}(t_1,\dots, t_k),\end{equation}
where 
\begin{equation}\label{continous}M_k(t_1,t_2,\cdots, t_k) = \int  \hat{\sigma}_{t_{k}}(\xi)  \widehat{f_{k-1}\mu}(-\xi)\hat{\mu}(\xi ) d \xi\end{equation}
is continuous and bounded below by a positive constant (independent of $\epsilon$) on $\tilde{I}\times \cdots \times \tilde{I}$, for a non-empty open interval $\tilde{I}$,
and
\begin{align}\label{remainder}
R_{k,j}^{\epsilon}(t_1,t_2,\cdots, t_k) 
=& \int \hat{\sigma}(t_{j}\xi) \left(1-\hat{\rho}(\epsilon\xi) \right) \widehat{f_{j-1}\mu}(\xi)\widehat{g^{\epsilon}_{j+1}\mu}(-\xi)d \xi\\
=& \mathcal{O}\left(\epsilon^{\alpha\left(s-\frac{d+1}{2}\right)}\right)
\end{align}
for some $\alpha>0$.
\end{lemma}

\vskip.125in 
In proving the lemma, we utilize the following notation:
\begin{equation}\label{back}g_j^{\epsilon}(x) = T_j^{\epsilon}\circ \cdots \circ T_k^{\epsilon}(1)(x),\end{equation} and $$g_{k+1}(x)=1.$$
It is important to note that $g_j(x)$ depends implicitly on  the choices of $t_1, \dots, t_k>0$, and this choice will be made explicit throughout.

First, we demonstrate equation \eqref{rewrite} with repeated use of Fourier inversion. 
We again employ a variant of the argument in \cite{IMT12}.
Write
\begin{align*}
\int f_{k}^\epsilon (x)d\mu(x) &= \int \int \sigma_{t_{1}}^\epsilon (x-y) g^{\epsilon}_{2}(y)d\mu(x)d\mu(y)\\
&=\int \int (\sigma_{t_{1}} * \rho_\epsilon)(x-y) g_{2}^{\epsilon}(y)d\mu(x)d\mu(y).
\end{align*}

\vskip.125in 
\noindent
Using Fourier inversion and properties of the Fourier transform, this is equal to

$$\int \int \int e^{2\pi i(x-y) \cdot \xi}\hat{\sigma}_{t_{1}}(\xi)\hat{\rho}_\epsilon(\xi) g_{2}^{\epsilon}(y)d\mu(x)d\mu(y)d \xi.$$

\vskip.125in 
\noindent
Simplifying further, we write
\begin{align*}
& \int f_{k}^\epsilon (x)d\mu(x)\\
=&\int \hat{\sigma}_{t_{1}}(\xi) \hat{\rho}(\epsilon\xi) \hat{\mu}(\xi)\widehat{g_{2}^{\epsilon}\mu}(-\xi)d \xi
\\
=&\left(\int \hat{\sigma}_{t_{1}}(\xi)\hat{\mu}(\xi)\widehat{g_{2}^{\epsilon}\mu}(-\xi)d \xi\right)
+ \left(\int \hat{\sigma}_{t_{1}}(\xi) \left( 1- \hat{\rho}(\epsilon\xi) \right) \hat{\mu}(\xi)\widehat{g_{2}^{\epsilon}\mu}(-\xi)d \xi\right)
\\
=&\left(\int \hat{\sigma}_{t_{1}}(\xi)\hat{\mu}(\xi)\widehat{g_{2}^{\epsilon}\mu}(-\xi)d \xi\right)
+ R_{k,1}^{\epsilon}(t_1,t_2,\cdots, t_k)
\end{align*}

With repeated use of Fourier inversion, we get
\begin{align*}
& \int f_{k}^\epsilon (x)d\mu(x)\\
=&\left( \int\hat{ \sigma}_{t_j}(\xi) \cdot \widehat{ f_{j-1}\mu}(-\xi)   
  \cdot \widehat{g^{\epsilon}_{j+1}\mu}(\xi)  d\xi\right)
+ \sum_{l=1}^j R_{k,l}^{\epsilon}(t_1,t_2,\cdots, t_k)\\
=& \, \, \cdots \\
=&\left( \int\widehat{ \sigma}_{t_k}(\xi)  \cdot  \widehat{f_{k-1}\mu}(-\xi)  
  \cdot \widehat{\mu}(\xi)     d\xi\right)
+ \sum_{l=1}^k R_{k,l}^{\epsilon}(t_1,t_2,\ldots, t_k)\\
=& \, \,M_k(t_1,t_2,\cdots, t_k)  \, \,
+ \, \, \sum_{l=1}^k R_{k,l}^{\epsilon}(t_1,t_2,\ldots, t_k) 
\end{align*}

We now prove that $M_k(t_1,t_2,\ldots, t_k)$ is continuous on any compact set away from $(t_1,\dots,t_k)=\vec{0}$ and that 
\begin{equation}\label{smallr}R_{k,j}^{\epsilon}(t_1,\ldots, t_k)= \mathcal{O}\left(\epsilon^{\alpha\left(s-\frac{d+1}{2}\right)}   \right).\end{equation}
Once these are established, we observe that the lower bound on constant chains established in \eqref{upperconstant} combined with \eqref{smallr} imply that $M_k(t_1,\dots,t_k)$ is positive when $t_1 =\cdots =t_k =t$ for any given $t \in I$. Fixing any such $t\in I$, it will then follow by continuity that $M_k(t_1,\dots, t_k)$ is bounded from below on $\tilde{I}\times \cdots \times \tilde{I}$ where $\tilde{I}$ is a non-empty interval.

\vspace{.125 in}

We now use the Dominated Convergence Theorem to verify the continuity of $M_k(t_1, \dots, t_k)$  on any compact set away from $(t_1,\dots,t_k)=\vec{0}$.
Let $t_1, \cdots, t_k>0$.  
Using properties of the Fourier transform and recalling the definition of $f_j$ from \eqref{fk} and $g_j$ from \eqref{back},  we write
$$M_k(t_1,t_2,\cdots, t_k)=\int\widehat{ \sigma}_{t_j}(\xi)  \cdot  \widehat{f_{j-1} \mu}(-\xi)     \cdot \widehat{ g_{j+1}\mu}(\xi)  d\xi$$
for any $j \in \left\{ 1, \ldots, k \right\}$.
\vskip.125in

Let $h_1, \ldots, h_k \in \mathbb{R}$ so that $(h_1, \ldots, h_k) \downarrow 0$.  
Let $$\tilde{f_j}= T_{t_{j} + h_{j}}\circ \cdots \circ T_{t_1 + h_1}(1) $$
and 
$$\tilde{g_j}=   T_{t_{j}+ h_{j}}\circ \cdots \circ T_{t_k + h_k}(1).$$

We have 
\\
\begin{align*} 
&M_k(t_1 + h_1,t_2 + h_2,\cdots, t_k + h_k)\\
&=\int\widehat{ \sigma}_{t_j + h_j}(\xi)  \cdot  \widehat{\tilde{f}_{j-1}\mu}(-\xi)     \cdot \widehat{\tilde{g}_{j+1} \mu}(\xi)  d\xi.
\\
\end{align*}
The integrand goes to $0$ as $h_j$ goes to $0$.   Now for $t_j$ in a compact set, the expression above is bounded by 
$$C(t_j) \int |\xi|^{-(d-1)/2} \left| \widehat{\tilde{f}_{j-1} \mu}(-\xi) \right|
\left| \widehat{\tilde{g}_{j+1} \mu}(\xi)\right|   d\xi.
$$
To proceed, we will utilize the following calculation. 

\begin{lemma} \label{fenergy} Let $\mu$ be a compactly supported Borel measure such that $\mu(B(x,r)) \leq Cr^s$ for some $s \in (0,d)$. Suppose that $\alpha>d-s$. Then for $f \in L^2(\mu)$, 
\begin{equation} \label{hibob} \int {|\widehat{f\mu}(\xi)|}^2 {|\xi|}^{-\alpha} d\xi \leq C'{||f||}^2_{L^2(\mu)}. \end{equation}
\end{lemma} 
\vskip.25in 

To prove Lemma \ref{fenergy}, observe that 
\begin{equation} \label{fenergysetup} \int {|\widehat{f\mu}(\xi)|}^2 {|\xi|}^{-\alpha} d\xi=C \int \int f(x)f(y) {|x-y|}^{-d+\alpha} 
d\mu(x)d\mu(y)=\langle Tf,f \rangle, \end{equation} where 
$$ Tf(x)=\int {|x-y|}^{-d+\alpha} f(y)d\mu(y)$$ and the inner product above is with respect to $L^2(\mu)$.
The positive constant, $C$, appearing in $\eqref{fenergysetup}$ depends only on the ambient dimension, $d$.  
 Observe that 
$$ \int {|x-y|}^{-d+\alpha} d\mu(y) \approx \sum_{j>0} 2^{j(d-\alpha)} \int_{|x-y| \approx 2^{-j}} d\mu(y) \leq C \sum_{j>0} 
2^{j(d-\alpha-s)} \leq C'$$ since $\alpha>d-s$. 

\vskip.125in 

By symmetry, $\int {|x-y|}^{-d+\alpha} d\mu(x) \leq C'$.  It follows by using Schur's test (\cite{Schur11}, see also Lemma 7.5 in \cite{W04}) that
$$ {||Tf||}_{L^2(\mu)} \leq C' {||f||}_{L^2(\mu)}.$$ 
\vskip.125in
This implies that conclusion of Lemma \ref{fenergy} by applying the Cauchy-Schwarz inequality to (\ref{fenergysetup}). This completes the proof of Lemma \ref{fenergy}. We note that Lemma \ref{fenergy} can also be recovered from the fractal Plancherel estimate due to R. Strichartz \cite{Str90}. See also Theorem 7.4 in \cite{W04} where a similar statement is proved by the same method as above. 

We already established using Theorem \cite{ISTU14} that finite compositions of the operators $T_l$ applied to $L^2(\mu)$ functions are in $L_2(\mu)$.  
Using the Cauchy-Schwarz inequality and in light of Lemma \ref{fenergy}, we $M_k(t_1 + h_1,t_2 + h_2,\cdots, t_k + h_k)$ is bounded.  We proceed by applying the Dominated Convergence Theorem.  
We have
\begin{align*}
&\lim_{h_j\downarrow 0} M_k(t_1 + h_1,t_2 + h_2,\cdots, t_k + h_k)\\
& =\int\widehat{ \sigma}_{t_j}(\xi)  \cdot  \widehat{\tilde{g}_{j-1} \mu}(-\xi)     \cdot
\widehat{ \tilde{f}_{j+1} \mu}(\xi)  d\xi.\\
& =\int\widehat{ \sigma}_{t_j}(\xi)  \cdot  \left( T_{t_{j-1} + h_{j-1}}\circ \cdots \circ T_{t_1 + h_1}(1)\cdot \mu\right)^{\widehat{}}(-\xi)     \cdot
\left(  T_{t_{j+1}+ h_{j+1}}\circ \cdots \circ T_{t_k + h_k}(1)\cdot \mu\right)^{\widehat{}}(\xi)  d\xi.\\
\end{align*}
We then rewrite the procedure, isolating $\widehat{\sigma}_{t_j}$ for each $j \in \{1, \dots, k\}$, and repeat the process above a total of $k$ times.  

\vskip.125in 
\noindent
\textbf{Bounding the remainder:}

Next, we wish to show that $\lim_{\epsilon \downarrow 0} R_k^\epsilon(t_1,\cdots, t_k) = 0$. Fix $\epsilon>0$.  Recall that $R_k^\epsilon(t_1,\cdots, t_k)$ is equal to
$$\int (1-\hat{\rho}(\epsilon \xi) )\hat{\sigma}(t\xi)\hat{\mu}(\xi)\widehat{f_k \mu}(-\xi)d \xi.$$
We consider the integral over $|\xi|<\left(\frac{1}{\epsilon}\right)^{\alpha}$ and the integral over $|\xi|>\left(\frac{1}{\epsilon}\right)^{\alpha}$ separately, where $\alpha \in (0,1)$ will be determined.  
Assume that $s>\frac{d+1}{2}.$

\begin{lemma}\label{lipschitz} Let $\rho:\mathbb{R}^d\rightarrow \mathbb{R}$ satisfy the following properties:  $\rho\geq 0$, $\rho(x) =\rho(-x)$, the support of $\rho$ is contained in $\{x: |x|<c\}$, and $\int \rho=1$.  Then
$$0\le 1-\widehat{\rho}(\xi) \le 2\pi c|\xi|.$$
\end{lemma}
To prove the Lemma \eqref{lipschitz}, write $$\widehat{\rho}(\xi) = \int \cos(2\pi x\cdot \xi) \rho(x) dx .$$ We observe that $\cos(x) + |x|>1$, and conclude that the lemma follows when $|x|<c.$  \\
It follows that
\begin{align*}
& \int_{|\xi|<\left(\frac{1}{\epsilon}\right)^{\alpha}}     \left|\hat{\rho}(\epsilon \xi) - 1\right| |\hat{\sigma}(t\xi)| |\hat{\mu}(\xi)||\widehat{f_k\mu}(-\xi)|d \xi   \\
&\lesssim \epsilon^{1-\alpha}\left(  \int  |\hat{\sigma}(t\xi)| |\hat{\mu}(\xi)||\widehat{f_k\mu}(-\xi)|d \xi  \right) \\
 & \lesssim\epsilon^{1-\alpha},\end{align*}
where the last line is justified in the estimation of $M_k(t)$ above.

\vspace{.25 in}
It remains to estimate the quantity

$$\int_{|\xi|>\left(\frac{1}{\epsilon}\right)^{\alpha}}   |\widehat{\sigma}(t\xi)| |\widehat{\mu}(\xi)| |\widehat{f_k\mu}(-\xi)| d\xi.$$
Proceeding as in the estimation of $M_k(t)$ above, we bound this integral above with
$$Ct^{-\frac{d-1}{2}}\int_{|\xi| >\left(\frac{1}{\epsilon}\right)^{\alpha}} |\xi|^{-\frac{d-1}{2}}|\hat{\mu}(\xi)||\widehat{f_k\mu}(-\xi)|d \xi$$
and then use Cauchy-Schwarz to bound it further with 
$$Ct^{-\frac{d-1}{2}}\left( \int_{|\xi| > \left(\frac{1}{\epsilon}\right)^{\alpha}} |\xi|^{-\frac{d-1}{2}}|\hat{\mu}(\xi)|^2 d \xi\right)^{1/2}\left( \int_{|\xi| > \left(\frac{1}{\epsilon}\right)^{\alpha}} |\xi|^{-\frac{d-1}{2}}|\widehat{f_k\mu}(\xi)|^2 d \xi\right)^{1/2}.$$

\vskip.125in 
\noindent
We have already shown that the second integral is finite. The first integral is bounded by 

$$\sum_{j >\alpha \log_2(1/\epsilon)} 2^{-j(\frac{d-1}{2})} \int_{2^j \leq |\xi| <2^{j+1}} |\hat{\mu}(\xi)|^2 d\xi.$$
We may choose a smooth cut-off function $\psi$ such that the inner integral is bounded by 

$$\int {|\widehat{\mu}(\xi)|}^2 \widehat{\psi}(2^{-j}\xi) d\xi.$$
By Fourier inversion, this integral is equal to 
$$ 2^{dj} \int \int \psi(2^j(x-y)) d\mu(x) d\mu(y) \leq C2^{j(d-s)}.$$ 

Returning to the sum, we now have the estimate

$$C\sum_{j >\alpha\log_2(1/\epsilon)} 2^{-j(\frac{d-1}{2})} \cdot 2^{j(d-s)} \leq C \sum_{j > \alpha\log_2(1/\epsilon)} 2^{j(\frac{d+1}{2}-s)}.$$
As long as $s > \frac{d+1}{2}$, this is $<< \epsilon^{\alpha(s-\frac{d+1}{2})}$. Thus $R_k^\epsilon(t_1,\dots, t_k)$ tends to 0 with $\epsilon$ as long as $\dim_\mathcal{H}(E) > \frac{d+1}{2}$. 

\vskip.125in 

In conclusion we have 
\begin{equation} \label{goodshit} \, \lim_{\epsilon \downarrow 0} \int \left(\prod_{j=1}^{k} \sigma_{t_j}^{\epsilon}(x^{i+1}-x^i) d\mu(x^i)\right) d\mu(x^{k+1})>c_k>0\end{equation} 
 for all $t_j\in \tilde{I}$. 

\vskip.25in

To complete the proof of Theorem \ref{main}, it remains to verify that $E$ contains a non-degenerate $k-$chain with prescribed gaps.  This is the topic of the next section.  

\section{Non-degeneracy}\label{degeneracy}

An important issue we have not yet addressed is that the chains we have found may be degenerate.   As an extreme example, consider the case where $t_i = 1 $ for all $i$. Then included in our chain count are chains which simply bounce back and forth between two different points.  We now take steps to insure that we can indeed find chains with distinct vertices.
\vskip.125in

We verified above that there exists a non-empty open interval $\tilde{I}$ so that 
$$\lim _{\epsilon\downarrow 0} \int \left(\prod_{j=1}^{k}  \sigma_{t_j}^{\epsilon}(x^{i+1}-x^i) d\mu(x^i)\right) d\mu(x^{k+1})$$
is bounded above and below for $t_1, \dots, t_k \in \tilde{I}$.
The upper bound appears in \eqref{up} and the lower bound appears in \eqref{goodshit}.
\vskip.125in

From here onward, we fix $ t_1, \dots, t_k \in \tilde{I}$ and set $\vec{t} = (t_1,\dots,  t_k)$.  
We now define a non-negative Borel measure on the set of $k-$chains with the gaps $\vec{t}$.
Let 
 $\Lambda_{\vec{t}}^k$ denote a non-negative Borel measure 
defined as follows
$$\Lambda^k_{\vec{t}}(A) =\lim _{\epsilon\downarrow 0} \int_A \left(\prod_{j=1}^{k}  \sigma_{t_j}^{\epsilon}(x^{i+1}-x^i) d\mu(x^i)\right) d\mu(x^{k+1}),$$
where $A\subset E\times \cdots \times E$, the $(k+1)-$fold product of the set $E$.
\vskip.125in

It follows that $\Lambda_{\vec{t}}^k$ is a finite measure which is not identically zero:
\begin{equation}\label{nonzed}0< \Lambda^k_{\vec{t}}(E\times \cdots \times E )        .\end{equation}
\vskip.125in

The strategy we use to demonstrate the existence of non-degenerate $k-$chains in $E$ is as follows: We first show that $\Lambda_{\vec{t}}^k$ has support contained in the set of $k-$chains.  This is accomplished by showing that the measure has support contained in all ``approximate'' $k-$chains.  We then show that the measure of the set of degenerate chains is zero.  It follows, since the $\Lambda_{\vec{t}}^k$-measure of the set of $k-$chains is positive and the $\Lambda_{\vec{t}}^k$-measure of the set of degenerate $k-$chains is zero, that the set of non-degenerate $k-$chains in $E$ is non-empty.  
\vskip.125in

For each $n \in \mathbb{Z}^+$, define the sets of $\frac{1}{n}-$approximate $k$-chains and the set of exact $k$-chains as follows:
$$A_{n,k}=\left\{\left(x^1, \dots, x^{k+1}\right)\in E\times \cdots \times E :t_i-\frac{1}{n}\le |x^{i+1} -x^i| \le t_i+ \frac{1}{n}, \text{ for each } i=1, \dots, k\right\},$$
and
$$A_k=\left\{\left(x^1, \dots, x^{k+1}\right)\in E\times \cdots \times E: |x^{i+1} -x^i| =t_i \text{ for each } i=1, \dots, k\right\}.$$
\vskip.125in

Observe that $$\bigcap_n A_{n,k}=A_k.$$
\vskip.125in

We now observe that the support of $\Lambda_t^k$ is contained in the set of all approximate chains.  
This follows immediately from the observation that  
$$\Lambda^k_{\vec{t}}(  A_{n,k}^c) =0,$$
for each $n\in \mathbb{Z}^+$,
where $  A_{n,k}^c$ denotes the compliment of the set $  A_{n,k}$ in $ E\times \cdots \times E$.  

Next, we observe that the support of  $\Lambda_{\vec{t}}^k$ is contained in the set of exact chains.  
Indeed, it follows from the previous equation that 
$$\Lambda^k_{\vec{t}}\left(  \bigcup_n A_{n,k}^c\right) \le \sum_n \Lambda^k_{\vec{t}}(  A_{n,k}^c)=0 .$$
Recalling \eqref{nonzed}, we conclude that 
\begin{equation} 0< \Lambda^k_{\vec{t}}(E\times \cdots \times E )   = 
\Lambda^k_{\vec{t}}\left(  \bigcup_n A_{n,k}^c\right) +  \Lambda^k_{\vec{t}}\left(  \bigcap_n A_{n,k}\right)  ,\end{equation}
and so
$$\Lambda^k_{\vec{t}} (A_{k})=\Lambda^k_{\vec{t}}\left(  \bigcap_n A_{n,k}\right)>0 .$$
Since $ t_1, \dots, t_k \in \tilde{I}$ were chosen arbitrarily, we have shown that 
$\Lambda^k_{\vec{t}} (A_{k})>0 $ whenever $\vec{t} = (t_1, \dots, t_k)$ and $t_i \in \tilde{I}$. 
\\

We now verify that the set of degenerate chains has $\Lambda^k_{\vec{t}}-$measure zero.  

\begin{lemma}
 Let
$$D_k = \{(x^1, . . ., x^{k+1})\in E\times \cdots \times E: x^i = x^j \text{ for some } i \neq j\}.$$
Then 
$$\Lambda^k_{\vec{t}}(D_k) =0.$$
\end{lemma}

\vspace{.25 in}

To prove the lemma, we first investigate the quantity
$$ \int_{D_k} \left(\prod_{j=1}^{k} \sigma_{t_j}^{\epsilon}(x^{i+1}-x^i) d\mu(x^i)\right) d\mu(x^{k+1}).$$
By the definition of $D_k$, we can bound this quantity above by

$$\sum_{1 \leq m <n\leq k+1} \int_{\{(x^1, \ldots , x^{k+1} : x^m = x^n\}} \left(\prod_{j=1}^{k} \sigma_{t_j}^{\epsilon}(x^{i+1}-x^i) d\mu(x^i)\right) d\mu(x^{k+1}).$$
We can rewrite the integral as

$$\int_{(\mathbb{R}^d)^k}  \int_{\{x: x = x^m\}} \left(\prod_{j=1}^{k} \sigma_{t_j}^{\epsilon}(x^{i+1}-x^i) \right) d\mu(x^n) d\mu(x^1) \cdots d\mu(x^{n-1})d\mu(x^{n+1})\cdots d\mu(x^{k+1}).$$
Since the inside integral is taken over a region of measure 0, this whole integral must be $0$. 
This holds for every choice of $m$ and $n$, and thus the entire sum must be $0$. This completes the proof of the lemma.
\\

In conclusion, we have shown that the set of exact $k-$chains has positive measure, $\Lambda^k_{\vec{t}}(A_k)>0$, and that the set of degenerate chains has zero measure,  $\Lambda^k_{\vec{t}}(D_k) =0$. It follows that $A_k\neq D_k$ and $A_k\neq \emptyset$.  In other words, there exists a non-empty open interval $\tilde{I}$, and there exists \textit{distinct} elements $x^1, \cdots, x^{k+1} \in E$ so that $|x^{i+1}-x^i|=t^i$ for each $i\in \{1, \dots, k\}$.
\\

\end{document}